\documentclass[12pt,a4paper,leqno]{amsart}
\usepackage{amsfonts,amssymb,amsmath,amsthm,graphicx}
\usepackage{a4wide}
\usepackage{mathrsfs}
\usepackage{epic,eepic}

\newtheorem{theorem}{Theorem}[section]
\newtheorem{proposition}[theorem]{Proposition}
\newtheorem{corollary}[theorem]{Corollary}
\newtheorem{lemma}[theorem]{Lemma}
\theoremstyle{remark}
\newtheorem{remark}[theorem]{Remark}
\numberwithin{equation}{section}

\newcommand{\N}{{\mathbb N}}
\newcommand{\R}{{\mathbb R}}
\newcommand{\E}{{\mathcal E}}
\newcommand{\D}{{\mathcal D}}
\renewcommand{\L}{\mathscr{L}}
\newcommand{\eps}{{\epsilon}}

\newcommand{\tx}[1]{\mbox{\;\;{#1}\;\;}}
\newcommand{\beqn}{\begin{eqnarray}}
\newcommand{\eeqn}{\end{eqnarray}}   
\newcommand{\beq}{\begin{eqnarray*}}
\newcommand{\eeq}{\end{eqnarray*}}
\newcommand{\hsub}{\underline{h}}
\newcommand{\Hsub}{\underline{H}}
\newcommand{\hsup}{\overline{h}}
\newcommand{\Hsup}{\overline{H}}
\newcommand{\usup}{\overline{u}}
\newcommand{\usub}{\underline{u}}

\begin{document}

\title[Large solutions to semilinear elliptic equations]{Large solutions to semilinear elliptic equations with Hardy potential and exponential nonlinearity}
\author{Catherine Bandle, Vitaly Moroz and Wolfgang Reichel}
\address{C. Bandle \hfill\break Mathematisches Institut, Universit\"at
  Basel \hfill\break
  Rheinsprung 21, CH-4051 Basel, Switzerland}
  \email{catherine.bandle@unibas.ch}
\address{V. Moroz \hfill\break
  Department of Mathematics, Swansea University \hfill\break
  Singleton Park, Swansea SA2 8PP, Wales, U.K.}
  \email{v.moroz@swansea.ac.uk}
\address{W. Reichel \hfill\break Fakult\"at \"fur Mathematik Universit\"at Karlsruhe, \hfill\break D-76128 Karlsruhe, Germany}
  \email{wolfgang.reichel@math.uni-karlsruhe.de}

\date{\today}

\begin{abstract}
On a bounded smooth domain $\Omega\subset \R^N$ we study solutions of a semilinear elliptic equation with an exponential nonlinearity and a Hardy potential depending on the distance to $\partial\Omega$. We derive global a priori bounds of the Keller--Osserman type. Using a Phragmen--Lindel\"of alternative for generalized sub and super-harmonic functions we discuss existence, nonexistence and uniqueness of so-called {\em large} solutions, i.e., solutions which tend to infinity at $\partial\Omega$. The approach develops the one used by the same authors \cite{BMR08} for a problem with a power nonlinearity instead of the exponential nonlinearity.
\end{abstract}


\maketitle

\section{Introduction}
Let $\Omega\subset \mathbb{R}^N$ be a bounded smooth domain (say $\partial\Omega\in C^3$) and let $\delta(x)$ be the distance from a point $x\in \Omega$ to the boundary $\partial \Omega$.  In this paper we study semilinear problems of the form
\begin{equation}\label{original}
-\Delta u - \frac{\mu}{\delta^2}u + e^u =0 \tx{in} \Omega,
\end{equation}
where $\mu\in\R$ is a given constant. The case without Hardy potential
\begin{align}\label{exponential}
-\Delta u + e^u = 0\tx{in} \Omega
\end{align}
is well-understood. In particular for any continuous function $\varphi\in C(\partial\Omega)$ the boundary value problem \eqref{exponential} with $u=\varphi \tx{on} \partial \Omega$
has a unique classical solution. Moreover there exists a unique  solution of \eqref{exponential}, cf. e.g. \cite{BM-JAM}, \cite{BM-CV}, with the property that
\begin{equation}\label{d-large}
u(x) \to \infty \tx{as} x\to \partial \Omega.
\end{equation}
This solution dominates all other solutions and is therefore commonly called {\sl large}. Near the boundary it behaves like \cite{BM-CV}
\begin{align} \label{bobe}
u(x)=\log \frac{2}{\delta^2(x)} +(N-1)\mathcal{H}_0(\sigma(x))\delta(x) +o(\delta(x)) \tx{as} x\to \partial \Omega,
\end{align}
where $\sigma:\Omega\to\partial\Omega$ denotes the nearest-point projection of $x$ onto the boundary and $\mathcal{H}_0(y)$ is the mean curvature of the boundary at $y\in\partial \Omega$.

\medskip

The presence of a Hardy potential has a significant effect on the set of solutions of \eqref{original}.
Because of the singularity of the potential the boundary values $\varphi$ in the problem
\begin{equation}\label{E-phi}
-\Delta u - \frac{\mu}{\delta^2}u + e^u =0 \mbox{ in } \Omega, \quad u=\varphi \mbox{ on } \partial \Omega
\end{equation}
cannot in general be prescribed arbitrarily. For instance, it is not difficult to show (see Theorem \ref{t-0} below)
that if $\varphi=0$ then problem \eqref{E-phi} admits a unique solution for every $\mu<C_H(\Omega)$, where
$C_H(\Omega)>0$ is the optimal constant in the Hardy's inequality
$$\int_\Omega|\nabla\phi|^2\,dx\ge C_H(\Omega)\int_\Omega\frac{\phi^2}{\delta^2}\,dx,\qquad\forall \phi\in C^\infty_0(\Omega).$$
On the other hand, if $\varphi>0$ is continuous then problem \eqref{E-phi} has no solution unless $\mu=0$.
This can be seen as follows. Without loss of generality let us assume that $u$ is positive in $\Omega$ (otherwise replace $\Omega$ by a neighbourhood of $\partial\Omega$). Suppose for contradiction that \eqref{E-phi} has a $C^2(\Omega)\cap C(\overline{\Omega})$-solution. Then the problem
\begin{equation}
-\Delta v = \frac{u}{\delta^2} \mbox{ in } \Omega, \quad v=0 \mbox{ on } \partial\Omega
\label{nonex}
\end{equation}
has a $C^2(\Omega)\cap C(\overline{\Omega})$-solution, where $v=\frac{1}{\mu}(u+z-h)$, $z$ is the Newtonian-potential of $e^u$ and $h$ is the harmonic extension of $(\varphi+z)|_{\partial\Omega}$. Let $f_k(x):=\min\{\frac{u(x)}{\delta^2(x)},k\}$ for $k\in \N$ and let $v_k$ be the weak $H_0^1(\Omega)$-solution of $-\Delta v_k = f_k$ in $\Omega$ with $v_k=0$ on $\partial\Omega$. Then $v_k\in C(\overline{\Omega})$ and
$$
v_k(x) = \int_\Omega G(x,y) f_k(y)\,dy \mbox{ for all } x \in \Omega,
$$
where $G(x,y)$ is the Dirichlet Green-function of $-\Delta$ on $\Omega$. The comparison principle yields $v_k(x)-\frac{1}{k} \leq u(x)$ for all $x\in \Omega$ and all $k\in \N$. However, by monotone convergence
$$
v_k(x) = \int_\Omega G(x,y) f_k(y)\,dy \to \int_\Omega G(x,y) \frac{u(y)}{\delta^2(y)} \,dy = \infty \mbox{ as } k \to \infty
$$
for all $x\in\Omega$. This is a contradiction.

\smallskip

The fact that no solutions exist with finite, non zero boundary data motivated us to study solutions which are unbounded near the boundary. The goal of the current paper is to study the {\sl large solutions} of \eqref{original}, i.e. solutions which satisfy \eqref{d-large}.
\smallskip

{\sc Main result.}
{\sl
i) If $\mu<0$ then \eqref{original} has no large solutions.
\smallskip

ii) If $0\leq \mu <C_H(\Omega)$ then there exists a unique large solution of \eqref{original}.
It is pointwise larger than any other solution of \eqref{original}}.
\bigskip

The paper is organized as follows.
In Section \ref{s-2} we set up the notation and introduce some basic definitions and tools. We also provide an existence
proof for the solution of \eqref{original}, vanishing on the boundary.
In Section \ref{s-3} we establish a Keller--Osserman type a priori upper bound on solutions of \eqref{original}.
In Section \ref{s-4} we prove the nonexistence of large solutions in the case $\mu<0$,
while in Sections \ref{s-5} and \ref{s-6} we establish asymptotic behavior, existence and uniqueness
of large solutions of \eqref{original} when $0\le\mu< C_H(\Omega)$.
Finally, in Section \ref{s-7} we construct a borderline case of a function $\gamma>0$ such that
$0<\gamma(\delta)\le 1$ and $\gamma(\delta)=o(\delta)$ as $\delta\to 0$ and for which the problem
$$
-\Delta u+\frac{\gamma(\delta)}{\delta^2}u+e^u =0 \tx{in} \Omega,
$$
has a large solution. We also discuss some open questions related to \eqref{original}.

\section{Some definitions and tools}\label{s-2}

For $\rho>0$ and $\eps\in(0,\rho)$ we use the notation
$$
\begin{array}{ll}
\Omega_\rho:=\{x\in\Omega:\delta(x)<\rho\},\quad & \Omega_{\eps,\rho}:=\{x\in\Omega:\eps<\delta(x)<\rho\}, \vspace{\jot}\\
D_\rho:=\{x\in \Omega: \delta(x)>\rho\},\quad
& \Gamma_\rho:=\{x\in\Omega:\delta(x)=\rho\}.
\end{array}
$$

\subsection{Sub- and super-harmonics}
For simplicity set
$$\L_\mu :=-\Delta - \frac{\mu}{\delta^2}.$$
Let $G\subset \Omega$ be open.
Following \cite{BMR08}, we call solutions $h$ of the equation
\begin{align}\label{harmonic}
\L_\mu h=0 \tx{in} G
\end{align}
{\sl harmonics} of $\L_\mu$ in $G$. If $G=\Omega$, we often omit $G$ and say that $h$ is a global harmonic of $\L_\mu$.
By interior regularity, weak solutions of \eqref{harmonic} are classical,
so in what follows we assume that all harmonics are of class $C^2(G)$.
\smallskip

We define {\sl super-harmonics} in $G$ as functions $\hsup\in H^1_{loc}(G)\cap C(G)$ which solve in the weak sense the differential inequality
\begin{equation}
\L_\mu \hsup \geq 0\tx{in}G.
\label{sub-super-harmonic}
\end{equation}
Similarly, $\hsub\in H^1_{loc}(G)\cap C(G)$ is called a {\sl sub-harmonic} in $G$ if the inequality sign is reversed.
\smallskip

If the functions $\hsub$ and $\hsup$ satisfy \eqref{sub-super-harmonic} in $\Omega$,
then they are called {\sl global sub or super-harmonics}, respectively.
If $\hsub$ and $\hsup$ satisfy \eqref{sub-super-harmonic} in a neighborhood of the boundary $\Omega_\eps$,
then they are respectively called {\sl local sub or super-harmonics}.
\smallskip

By the classical strong maximum principle for the Laplacian with potentials applied locally in small subdomains of $\Omega$,
any nontrivial super-har\-monic $\hsup \gneq 0$ is strictly positive in $\Omega$, while any sub-harmonic $\hsub$  in $\Omega$ is locally bounded above.
\smallskip

The following examples of explicit local sub and super-harmonics  will play an important role in our considerations.

\smallskip

\noindent
{\sc Examples} \cite[Lemma 2.8]{BMR08}. Let  $\mu<1/4$ and
$$
\beta_\pm = \frac{1}{2} \pm \sqrt{\frac{1}{4}-\mu}.
$$
The function $\delta^\beta$ is a local super-harmonic of $\L_\mu$ if $\beta\in (\beta_-,\beta_+)$. It is a local sub-harmonic if $\beta \notin [\beta_-,\beta_+]$. In the borderline cases $\beta=\beta_\pm$, we have that
for small $\eps>0$
$$\hsup= \delta^{\beta_+}(1-\delta^\eps), \quad \Hsup=\delta^{\beta_-}(1+\delta^\eps)
$$
are local super-harmonics and
$$\hsub= \delta^{\beta_+}(1+\delta^\eps), \quad \Hsub=\delta^{\beta_-}(1-\delta^\eps)
$$
are local sub-harmonics.

\subsection{Hardy constant}\label{C-Hardy}

The constant
$$
C_H(\Omega)= \inf_{0\neq \phi\in W^{1,2}_0(\Omega)} \frac{\int_\Omega |\nabla \phi|^2 \:dx}{\int_\Omega \delta^{-2}(x)\phi^2\:dx}
$$
is called the {\sl global Hardy constant}. It is well-known that $0<C_H(\Omega)\le 1/4$. In general $C_H(\Omega)$ varies with the domain. For convex domains $C_H(\Omega)=1/4$,
but there exist smooth domains for which $C_H(\Omega)<1/4$. A review with an extensive bibliography and where, in particular, Maz'ya's relevant earlier contributions \cite{Mazja} are mentioned, is found in \cite{Da99}. Improvements of this inequality by adding an additional $L^q$ norm were obtained
by Filippas, Maz'ya and Tertikas in a series of papers.
The most recent results are found in \cite{FiMaTe07}.
This paper contains also references to previous related works.
It turns out,  cf.~\cite{MMP}, that  $C_H(\Omega)$ is attained if and only if $C_H(\Omega)<1/4$.
Notice that $C_H(\Omega)$ is in general not monotone with respect to $\Omega$. 

The relation between Hardy's constant, existence of positive super-har\-monics in $\Omega$, and validity of a comparison principle
for $\L_\mu$ is explained by the following classical result (cf. \cite[Theorem 3.3]{Agmon}).
\begin{lemma}\label{t-GST}
The following three statements are equivalent:
\begin{itemize}
\item[$(i)$] $\mu\le C_H(\Omega)$.
\item[$(ii)$] $\L_\mu$ admits a positive super-harmonic in $\Omega$.
\item[$(iii)$] For any subdomain $G$ with $\overline{G}\subset\Omega$ and any two sub and super-harmonics $\hsub$, $\hsup$ of $\L_\mu$ in $G$ with $\hsub\leq \hsup$ on $\partial G$ it follows that $\hsub \leq \hsup$ a.e. in $G$.

\end{itemize}
\end{lemma}

\subsection{Phragmen--Lindel\"of alternative}
Observe that global positive super-harmonics of $\L_\mu$ exist for all $\mu\le C_H(\Omega)$,
while the existence of local positive super-harmonics of $\L_\mu$ is controlled by the {\sl local Hardy constant}
$$C_H^{loc}(\Omega_\rho):=
\inf_{0\neq \phi\in W^{1,2}_0(\Omega_\rho)} \frac{\int_{\Omega_\rho} |\nabla \phi|^2 \:dx}{\int_{\Omega_\rho} \delta^{-2}(x)\phi^2\:dx}.$$
Note that in generally, $C_H(\Omega_\rho)\neq C_H^{loc}(\Omega_\rho)$ because
$\delta(x)=\mathrm{dist}(x,\partial\Omega)\not = \mathrm{dist}(x,\partial \Omega_\rho)$.
It is known \cite[Lemma 2.5]{BMR08} that $C_H^{loc}(\Omega_\rho)=1/4$ if $\rho>0$ is sufficiently small. 

If $\mu\leq C_H^{loc}(\Omega_\rho)$ then $\L_\mu$ admits positive local super--harmonics and satisfies the comparison principle
between sub and super--harmonics in $\Omega_\rho$, for all sufficiently small $\rho>0$, see \cite{BMR08}.
Furthermore, the following {\sl Phragmen--Lindel\"of alternative} holds for $\L_\mu$.
We repeat the statement and its proof from \cite[Theorem 2.6]{BMR08} for the reader's convenience.


\begin{theorem}\label{t-PL}
Let $\mu\le 1/4$. Let $\hsub$ be a local positive sub-harmonic. Then the following alternative holds:
\begin{enumerate}
\item[{$(i)$}]
either for every local super-harmonic $\hsup>0$
\begin{equation}\label{i}
\limsup_{x\to\partial\Omega}\hsub/\hsup>0,
\end{equation}
\item[{$(ii)$}]
or for every local super-harmonic $\hsup>0$
\begin{equation}\label{ii}
\limsup_{x\to \partial\Omega} \hsub/\hsup<\infty.
\end{equation}
\end{enumerate}
\end{theorem}

\proof
Assume $(i)$ does not hold, that is there exists a super-harmonic $\hsup_\ast>0$ that
\begin{equation}\label{e-A-small}
\lim_{x\to\partial\Omega} \hsub/\hsup_\ast=0.
\end{equation}
Let $\hsup>0$ be an arbitrary super-harmonic in $\Omega_\rho$ for some sufficiently small $\rho>0$. Then there exists a constant $c>0$ such that $\hsup \ge c\hsub$ on $\Gamma_{\rho/2}$.
For $\tau>0$, define a comparison function
$$
v_\tau:=c\hsub-\tau \hsup_\ast.$$
Then \eqref{e-A-small} implies that for every $\tau>0$ there exists $\eps=\eps(\tau)\in(0,\rho)$ such that
$v_\tau\le 0$ on $\Omega_{\eps}$.
Applying the comparison principle in $\Omega_{\eps/2,\rho/2}$, we conclude that
$\hsup\ge v_\tau$ in $\Omega_{\eps/2,\rho/2}$ and hence, in $\Omega_{\rho/2}$.
So by considering arbitrary small $\tau>0$, we conclude that
for every super-harmonic $\hsup>0$ in $\Omega_\rho$ there exist $c>0$
such that $\hsup\geq c\hsub$ holds in $\Omega_\rho$. This implies \eqref{ii}.
\qed

\medskip

If we apply this alternative to the special super-harmonics mentioned above we get for sub-harmonics the following  boundary behavior. If $\mu <1/4$ then either
$$
(i) \quad \quad \quad \quad \limsup_{x\to \partial \Omega} \frac{\hsub(x)}{\delta(x)^{\beta_-}}>0
$$
or
$$
(ii) \quad \quad \quad \quad\limsup_{x\to \partial \Omega} \frac{\hsub(x)}{\delta(x)^{\beta_+}}<\infty.
$$


\subsection{Sub- and super-solutions}
Let $G\subset \Omega$ be open. A function $\usup\in H^1_{loc}(G)\cap C(G)$ satisfying the inequality
$$
\L_\mu \usup +e^{\usup} \geq 0 \tx{in} G
$$
in the weak sense is called a {\sl super-solution} of \eqref{original} on $G$. Similarly $\usub\in H^1_{loc}(G)\cap C(G)$ is called a {\sl sub-solution} of \eqref{original} if the inequality sign is reversed. A function $u$ is a solution of \eqref{original} in $G$ if it is a sub and super-solution in $G$.
By interior elliptic regularity weak solutions of \eqref{original} are classical.
Hence in what follows we assume that all solutions of \eqref{original} are of class $C^2(\Omega)$.
\medskip

Observe that solutions and sub-solutions are sub-harmonics of $\L_\mu$.

\medskip

The following comparison principle is based on an argument used in \cite{BMR08} and plays a crucial role in our estimates.
Part (i) relies heavily on the fact that $\mu<C_H(\Omega)$. Part (ii) is an extension of (i) for arbitrary $\mu$ under an additional assumption.

\begin{lemma}[\sc Comparison principle]
\label{comp2}
Let $G\subset \Omega$ be open and let $\usub, \usup \in H^1_{loc}(G)\cap C(G)$ be a pair of sub-, super-solutions to \eqref{original} satisfying
$$
\limsup_{x\to\partial G} [\usub(x)-\usup(x)] <0.
$$
\begin{itemize}
\item[(i)] If $\mu<C_H(\Omega)$ then $\usub\leq \usup$ in $G$.
\item[$(ii)$] If $\mu \geq C_H(\Omega)$ and in addition $\usup>1$ in $G$ then $\usub\leq \usup$ in $G$.
\end{itemize}
\end{lemma}
\proof Let $G_+ := \{x\in G: \usub(x)>\usup(x)\}$. In view of the boundary conditions we have  $\overline{G}_+\subset G$. In the weak formulation of the inequality
\begin{equation}
\L_\mu (\usup-\usub) \geq -(e^{\usup}-e^{\usub}) \tx{in} G
\label{weak_ineq}
\end{equation}
we use the test function $(\usub-\usup)_+\in H_0^1(G)$ and obtain
$$
\int_G |\nabla (\usub-\usup)_+|^2\,dx \leq \mu \int_G \delta^{-2}(\usub-\usup)_+^2 \,dx.
$$
Case (i): unless $G_+=\emptyset$ this implies
$$
\mu \geq \inf_{0\neq \phi\in W^{1,2}_0(G)} \frac{\int_G |\nabla \phi|^2 \:dx}{\int_G \delta^{-2}(x)\phi^2\:dx}\geq
\inf_{0\neq \phi\in W^{1,2}_0(\Omega)} \frac{\int_\Omega |\nabla \phi|^2 \:dx}{\int_\Omega \delta^{-2}(x)\phi^2\:dx}=
C_H(\Omega),
$$
which contradicts our assumption.

\smallskip

\noindent
Case (ii): if $\mu\geq C_H(\Omega)$ we make use of the following argument. In the weak formulation \eqref{weak_ineq} we use again the test function $(\usub-\usup)_+\in H_0^1(G)$ and obtain
\begin{align} \label{a1}
\int_G |\nabla (\usub-\usup)_+|^2\,dx -\mu\int_G \delta^{-2}(\usub-\usup)_+^2\, dx \leq \int_{G_+} \frac{e^{\usup} -e^{\usub}}{\usub-\usup}(\usub-\usup)^2_+ \,dx.
\end{align}
Since $\usup>1$ in $G$ we can write  $(\usub-\usup)_+=\phi \usup$ where $\phi\in W^{1,2}_0(G)$ and the support of $\phi$ lies in the closure of $G_+$. Then
\begin{align*}
\int_G |\nabla (\usub-\usup)_+|^2\, dx &= \int_{G_+} (\phi^2|\nabla \usup|^2 +2\phi\usup \nabla \usup\cdot\nabla \phi +\usup^2|\nabla \phi|^2)\,dx\\
&=\int_{G_+}[\usup^2|\nabla \phi|^2+ \nabla \usup\cdot\nabla(\phi^2\usup)]\,dx.
\end{align*}
Recalling that $\usup$ is a super solution and that $\phi^2\usup\geq 0$, we conclude that
$$
\int_{G_+} \nabla \usup \cdot \nabla(\phi^2\usup)\,dx\geq \int_{G_+}\big[\frac{\mu}{\delta^2}\phi^2\usup^2-e^{\usup}\phi ^2\usup\big]\,dx.
$$
This leads to
\begin{align}\label{a2}
\int_G |\nabla (\usub-\usup)_+|^2\, dx -\mu\int_{G_+}\delta^{-2}(\usub-\usup)^2 \,dx\geq -\int_{G_+}\frac{e^{\usup}}{\usup}(\usub-\usup)^2\,dx.
\end{align}
Since by convexity
$$
\frac{e^{\usup}-e^{\usub}}{\usub-\usup}\leq -e^{\usup} \tx{whenever} \usub\geq \usup,
$$
and moreover $\usup>1$ by assumption, we find that
\eqref{a2} contradicts \eqref{a1} unless $G_+=\emptyset$.
\qed

\subsection{Solutions with zero boundary data}
We are going to show that the problem
\begin{equation}\label{E-0}
\L_\mu u+e^u=0,\qquad u\in H^1_0(\Omega),
\end{equation}
admits a solution for all $\mu<C_H(\Omega)$. For this purpose we need the following lemma.
\begin{lemma}
\label{variational}
Let $\mu<C_H(\Omega)$. Then the boundary value problem
\begin{equation}
\L_\mu \phi=-1, \qquad \phi\in H^1_0(\Omega),
\label{torsion_like}
\end{equation}
admits a unique solution $\phi<0$. In addition  $\phi$ is bounded in $\overline{\Omega}$.
\end{lemma}

\proof
Results of this type are standard, cf. for instance \cite{MMP} and the references given there.
For the sake of completeness we sketch the proof.
Consider the quadratic form associated to $\L_\mu$:
$$\E_\mu(u):=\int_\Omega\Big(|\nabla u|^2
-\mu\frac{u^2}{\delta^2}\Big)dx.$$
It follows from the definition of the Hardy-constant $C_H(\Omega)$ that
\begin{equation}\label{coercive}
\E_\mu(u)\ge\Big(1-\frac{\mu}{C_H(\Omega)}\Big)\int_\Omega|\nabla u|^2 dx.
\end{equation}
We conclude that $\E_\mu$ is a coercive and continuous quadratic form on
$H^1_0(\Omega)$.
Since $-1\in [H^1_0(\Omega)]^\ast$, the existence and uniqueness of the
solution $\phi\in H^1_0(\Omega)$ follows by the Lax--Milgram theorem.
Since $\L_\mu \phi<0$ in $\Omega$, the comparison principle of Lemma \ref{t-GST}
implies that $\phi<0$. By the classical regularity theory $\phi$ is bounded in every compact subset of $\Omega$. A straightforward computation (using formula \eqref{mean}) shows that for large $A$ and small $\epsilon$
$$
\underline{\phi}=-A\delta^\nu, \quad \nu=\min\{2,\beta_+ -\epsilon\}
$$
is a sub-solution in $\Omega_{\delta_0}$ for a small $\delta_0>0$. By chosing $A>0$ so large that in particular $\underline{\phi} \leq \phi$ on $\Gamma_{\delta_0}$ we can apply the comparison principle and conclude that $\phi$ is bounded in $\overline{\Omega}$.
\qed

\begin{theorem}\label{t-0}
Let $\mu<C_H(\Omega)$. Then \eqref{E-0} has a unique solution $u_0$.
Moreover, $\phi<u_0<0$, where $\phi$ is defined in Lemma \ref{variational}.
\end{theorem}

\proof
Consider the  energy functional corresponding to \eqref{E-0}:
$$
J(u):=\frac{1}{2}\E_\mu(u)+\int_\Omega e^u\,dx.
$$
In view of \eqref{coercive}, it is clear that $J:H^1_0(\Omega)\to\R\cup\{+\infty\}$
is coercive, convex and weakly lower semicontinuous on $H^1_0(\Omega)$.
Hence $J$ admits the unique minimizer $u_0\in H^1_0(\Omega)$.
Note that $J(u^-)\le J(u)$ for every $u\in H^1_0(\Omega)$.
As a consequence, $u_0\le 0$. Hence $e^{u_0}$ is bounded from above, and thus $u_0$ satisfies the Euler-Lagrange equation and solves \eqref{E-0}.
Further, since $u=0$ is not a solution of \eqref{E-0} we conclude that $u_0<0$.

Let $\phi\in H^1_0(\Omega)$ be as defined in Lemma \ref{variational}.
Since $\phi<0$ in $\Omega$, we have $\L_\mu \phi+e^\phi\leq 0$ in $\Omega$,
so $\phi$ is a sub-solution of \eqref{E-0}.
>From the comparison principle of Lemma \ref{comp2} (i) it follows that $\phi<u_0$.
\qed

\begin{remark}
Suppose the domain $\Omega$ is such that $C_H(\Omega)<1/4$. Then there exists a positive solution $\phi_1\in H^1_0(\Omega)$ of $\L_{C_H(\Omega)} \phi_1=0$ in $\Omega$, see \cite{MMP}. We claim that if $\mu\geq C_H(\Omega)$ then \eqref{E-0} has no negative solution. Suppose $u\in H_0^1(\Omega)$ is a negative solution of \eqref{E-0}. Then we obtain the contradiction that
$$
0\leq \int_\Omega \frac{(C_{H(\Omega)}-\mu)}{\delta^2} u\phi_1\,dx = \int_\Omega\nabla u \cdot \nabla \phi_1- \frac{\mu}{\delta^2}u\phi_1\,dx= -\int_\Omega e^u\phi_1\:dx<0.
$$
Hence, if for $C_H(\Omega)<1/4$ and $\mu\geq C_H(\Omega)$ a solution of \eqref{E-0} exists then it must be sign-changing, cf. Question 1 in Section~\ref{s-7}. The same statement holds for solutions of \eqref{torsion_like}.
\label{rem}
\end{remark}


\section{A priori upper bounds}\label{s-3}
In  this section we construct a universal upper bound for all solutions of \eqref{original} by means of a super-solution which tends to infinity at the boundary. The construction is inspired by the Keller--Osserman bound given in \cite{BMR08} for power nonlinearities. The terminology {\em Keller-Osserman bound} refers the universal upper bound of Lemma~\ref{Ob1} and Lemma~\ref{Osserman}. Such upper bounds, which hold for all solutions of a nonlinear equation, were observed in the classical papers by Keller \cite{Keller} and Osserman \cite{Osser}.

\medskip

For our purpose we need the {\sl Whitney distance} $d:\Omega \to \mathbb{R}_+$ which is a $C^\infty(\Omega)$-function such that for all $x\in \Omega$
\begin{align*}
&c^{-1}\delta(x)\leq d(x)\leq c\delta(x),\\
&|\nabla d(x)|\leq c,\quad |\Delta d(x)|\leq cd^{-1}(x),
\end{align*}
with a constant $c>0$ which is independent of $x$. These properties of the Whitney distance may be found in \cite{S}.

For $\eps>0$, we use the notation $\D_\eps=\{x\in\Omega:d(x)>\eps\}$.
\begin{lemma}\label{Ob1}
Let $\mu\leq 0$. Then there exists a number $A>0$ such that for every solution $u$ of \eqref{original} we have
$$
u(x)\leq \log \frac{A}{d^2(x)} \mbox{ in } \Omega.
$$
\end{lemma}
\proof
Consider for small $\eps>0$ the function $f_\eps(x)= \log \frac{A}{(d(x)-\eps)^2}$ in $\D_\eps$. It satisfies
$$
\Delta f_\eps = \frac{2}{(d-\eps)^2} |\nabla d|^2 -\frac{2}{d-\eps}\Delta d \quad \mbox{ in } \D_\eps.
$$
Thus by the properties of the Whitney distance and since $\mu$ is non-positive
$$
\Delta f_\eps + \frac{\mu}{\delta^2}f_\eps -e^{f_\eps}\leq \frac{c_1-A}{(d-\eps)^2} \quad \mbox{ in } \D_\eps.
$$
For sufficiently large $A$, the right-hand side of this inequality is negative. Hence $f_\eps$ is a super-solution satisfying $f_\eps>u$ on $\partial \D_\eps$. The comparison principle implies that
$$u(x)\leq f_\eps(x)\tx{in} \D_\eps.$$
Since $\eps>0$ is arbitrary the conclusion follows.
\qed

\begin{figure}[h!t]
\centering
\includegraphics[width=6cm]{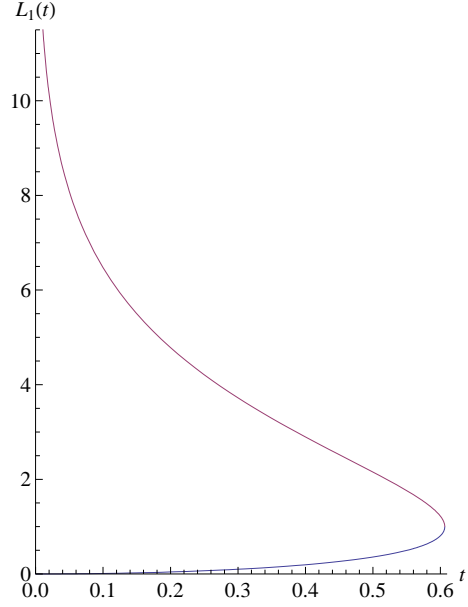}
\caption{Lambert function $L_1(t)$.}
\label{lambert_fig}
\end{figure}

If $\mu$ is positive we proceed differently.
For $A>0$ the function $L_A(d(x))$, $d(x)=$Whitney-distance, will play an essential role in the following construction of upper bounds for all solutions of \eqref{original}. The definition of $L_A(t)$ is given implicitly by
\begin{align}\label{lambert}
\frac{e^{L_A(t)}}{L_A(t)}= \frac{A}{t^2}, \quad A>0, \quad L_A(t)>1.
\end{align}
It is easily seen that $L_A(t)$ is defined whenever $A\geq et^2$ and that it has two branches. We select the branch $L_A(t)\geq 1$, cf. Figure \ref{lambert_fig}.
Clearly the function $L_A(t)$ is monotone increasing in $A$ and decreasing in $t$. Also, from the relation $L_A(t) = \log\frac{A}{t^2}+\log L_A(t)$ one finds successively
\begin{align*}
L_A(t) & \geq \log\frac{A}{t^2}, \\
L_A(t) & \geq \log\frac{A}{t^2}+\log\log\frac{A}{t^2}, \\
L_A(t) & \geq \log\frac{A}{t^2}+\log\left(\log\frac{A}{t^2}+\log\log\frac{A}{t^2}\right), \\
L_A(t) & \geq \ldots
\end{align*}
Moreover
\begin{align}\label{L_as}
 \lim_{t\to 0+} \frac{L_A(t)}{\log(1/t^2)}=\lim_{t\to 0+}
\frac{L_A'(t)}{-2/t}=\lim_{t\to 0+} \frac{L_A(t)}{L_A(t)-1}=1,
\end{align}
since $L_A(0)=\infty$.

As a historical note, let us mention that the function $L_A(t)$ is related to the Lambert $W$-function which satisfies the equation
$$
W(s)e^{W(s)}=s
$$
and which has a long history starting with J.H.~Lambert and L.~Euler. Indeed we have
$$
L_A(t) = -W(-\frac{t^2}{A}),
$$
if one takes for $W$ again the upper branch.

\medskip

Next we show that $L_A(d(x))$ is indeed a universal upper bound for all solutions of \eqref{original} provided one takes $A>0$ sufficiently large. The estimate is based on the extended comparison principle of Lemma~\ref{comp2}(ii).

\begin{lemma}\label{Osserman}
There exists $A>0$ such that every solution of \eqref{original} satisfies
$$
u(x)\leq L_A(d(x)) \tx{in} \Omega.
$$
\end{lemma}
\proof
In order to define $L_A(d(x))$ with property \eqref{lambert} we must take $A>0$ so large that $\inf_{\Omega}\frac{A}{d^2(x)} >e$. A straightforward computation yields
\begin{align}\label{Laplacian}
\Delta L_A(d)= \frac{2L_A(d)}{d^2(L_A(d)-1)}|\nabla d|^2\left\{1-\frac{2}{(L_A(d)-1)^2}\right\}-\frac{2L_A(d)}{d(L_A(d)-1)}\Delta d.
\end{align}
For $\eps\ge 0$, let $\usup_\eps: \D_\eps\to \R$ be defined as
\begin{align}\label{e-usup-eps}
\usup_\eps(x):=L_A(d(x)-\eps).
\end{align}
Then by \eqref{lambert}, \eqref{Laplacian} and the properties of the Whitney distance
\begin{align*}
\Delta \usup_\eps +\frac{\mu}{\delta^2}\usup_\eps -e^{\usup_\eps} \leq & \frac{2L_A(d-\eps)}{(d-\eps)^2(L_A(d-\eps)-1)}c^2\left\{1-\frac{2}{(L_A(d-\eps)-1)^2}+\frac{2}{c}\right\} \\
& +\frac{L_A(d-\eps) }{(d-\eps)^2}\left\{c_\pm \mu-A\right\},
\tx{where} c_\pm=
\begin{cases}
0 &\tx{if} \mu\leq 0\\
c^2&\tx{if} \mu>0.
\end{cases}
\end{align*}
By taking $A$ sufficiently large we can always achieve that the right-hand side is negative independently of $\eps$. Consequently $\usup_\eps$ is a super-solution of \eqref{original} in $\D_\eps$, for all sufficiently small $\eps> 0$.

Let $u$ be an arbitrary solution of \eqref{original}.
Clearly $u<\usup_\eps$ on $\partial \D_\eps$.
Moreover, by definition $\usup_\eps(x)=L_A(d(x)-\eps)>1$ and thus Lemma~\ref{comp2}(ii) applies and yields
$$
u(x)\leq L_A(d(x)-\eps) \tx{in} \D_\eps.
$$
Since $\eps>0$ was an arbitrary small number, this concludes the proof of the lemma.
\qed

\begin{remark}
It is clear from the above proof that for a sufficiently large $A>0$ the function
\begin{align}\label{e-usup}
\usup(x):=L_A(d(x))
\end{align}
is a super-solution of equation \eqref{original} in $\Omega$.
\end{remark}

\begin{remark}
Notice that
$$
 \frac{e^{L_A(d(x))}}{L_A(d(x))}=\frac{A}{d^2(x)}\leq \frac{Ac^2}{\delta^2(x)}= \frac{e^{L_{c^2A}(\delta(x))}}{L_{c^2A}(\delta(x))}.
$$
Replacing the Whitney distance by the standard distance we obtain the universal a priori bound
$$
u(x) \leq L_{c^2A}(\delta(x)),
$$
and by \eqref{L_as} we obtain
\begin{align}\label{uAs}
\limsup_{x\to\partial\Omega}\frac{u(x)}{\log \delta^{-2}(x)}\leq 1.
\end{align}
It should be pointed out that the bound constructed above holds for every $\mu\in \mathbb{R}$.
\end{remark}

\section{Nonexistence of large solutions if $\mu<0$}\label{s-4}
Lemma \ref{Osserman} together with the Phragmen--Lindel\"of alternative gives rise to a nonexistence result.
\begin{theorem} \label{nonexistence} If $\mu<0$ then \eqref{original} does not have large solutions.
\end{theorem}
\proof
If a solution $u$ of \eqref{original} exists with $u(x)\to \infty$ as $x\to \partial \Omega$, then by conclusion drawn from the Phragmen--Lindel\"of alternative of Theorem \ref{t-PL} it must satisfy
$$
\limsup_{x\to \partial \Omega} \frac{u(x)}{\delta(x)^{\beta_-}}>0, \tx{where} \beta_-= \frac{1}{2} - \sqrt{\frac{1}{4}-\mu}.
$$
On the other hand \eqref{uAs} implies
$$
\limsup_{x\to \partial \Omega} \frac{u(x)}{\delta(x)^{\beta_-}}\leq \limsup_{x\to \partial \Omega}\delta(x)^{-\beta_-}\log \frac{1}{\delta(x)^2} =0.
$$
This is impossible and therefore $u$ does not exist.
\qed

This nonexistence result together with the Phragmen--Lindel\"of alternative leads to the following conclusion.
\begin{corollary}
If $\mu<0$ then all solutions of \eqref{original} vanish on the boundary.
\end{corollary}
\section{Asymptotic behavior of large solutions near the boundary}\label{s-5}
\subsection{Global sub solutions}
Since the  case $\mu=0$ is well-known and since no large solutions exist for negative $\mu$ we shall assume throughout this section that $\mu>0$.

Let $L_A$ be defined as in \eqref{lambert}. Next we shall construct local sub-solutions which have the same asymptotic behavior as the super-solution $L_A(d(x))$ from Lemma~\ref{Osserman}.
\begin{proposition}
\label{local_sub}
Let $0<B\leq \mu$. Then there exists a small positive $\eps_0<\frac{1}{2}\sqrt{B/e}$ such that $\usub_\eps(x):=L_B(\delta(x)+\eps)$ is a sub solution of \eqref{original} in $\Omega_{\eps_0}$ for any $\eps\in [0,\eps_0]$.
\end{proposition}
\proof
Since $(\delta(x)+\eps)^2 \leq 4\epsilon_0^2<B/e$ the function $\usub_\eps$ is well defined in $\Omega_{\eps_0}$. We have, as in the proof of Lemma~\ref{Osserman}
\begin{align}\label{lower}
\Delta \usub_\eps +\frac{\mu}{\delta^2}\usub_\eps-e^{\usub_\eps}=\frac{2\usub_\eps}{(\delta+\eps)^2(\usub_\eps-1)}\left\{1-\frac{2}{(\usub_\eps-1)^2}\right\}-\frac{2\usub_\eps}{(\delta+\eps)(\usub_\eps-1)}\Delta \delta\\
\nonumber +\frac{\mu}{\delta^2}\usub_\eps-\frac{B}{(\delta+\eps)^2}\usub_\eps.
\end{align}
In $\Omega_{\eps_0}$ one has the expansion
\begin{equation}
\Delta \delta(x)= -(N-1){\mathcal H}_0(\sigma(x)) +o(\delta(x)).
\label{mean}
\end{equation}
and hence $\Delta\delta \leq K$ in $\Omega_{\eps_0}$ for some constant $K>0$ independently of $\eps_0$. Next we choose $\eps_0$ so small that
$1-\frac{2}{(\usub_\eps-1)^2}\geq \frac{1}{2}$ in $\Omega_{\eps_0}$. Since $0<B\leq\mu$ we find
$$
\Delta \usub_\eps +\frac{\mu}{\delta^2}\usub_\eps-e^{\usub_\eps} \geq
\frac{\usub_\eps}{(\delta+\eps)(\usub_\eps-1)}\left(\frac{1}{\delta+\eps}-2K\right)
$$
in $\Omega_{\eps_0}$. The right-hand side is positive provided $\epsilon_0<1/(4K)$. Thus $\usub_\eps$ is a sub-solution in $\Omega_{\eps_0}$ for all $\eps\in [0,\eps_0]$.
\qed

\medskip

In the next step we extend the local sub-solution $\usub_\eps$ to a global sub-solution $\underline{\mathcal U}_\eps$ in the whole domain such that
$\underline{\mathcal{U}}_\eps=\usub_\eps$ near the boundary.

\begin{proposition}\label{lbound}
Assume  $0<\mu<C_H(\Omega)$. Then there exists a global sub-solution $\underline{\mathcal U}_0$ with
$\underline{\mathcal U}_0(x)=L_\mu(\delta(x))\bigl(1-O(\delta(x)^{\beta_-})\bigr)$. Moreover, if $u$ is any solution of \eqref{original} which tends to infinity at the boundary then
$u\geq \underline{\mathcal{U}}_0$ and in particular
\begin{equation}
\label{limit}
\liminf_{x\to \partial \Omega} \frac{u(x)}{\log \delta^{-2}(x)}\geq 1.
\end{equation}
\end{proposition}
\proof
Let $\phi\in H^1_0(\Omega)$ be as defined in Lemma \ref{variational}. Since $\phi$ is non positive, we have $\L_\mu \phi+e^\phi\leq 0$ in $\Omega$ and $\phi$ is therefore a sub-solution of \eqref{original}. Let $\usub_\eps(x) = L_\mu(\delta(x)+\eps)$ by the local sub-solution from Proposition \ref{local_sub}. Consider the local super-harmonic (cf. {\sc Examples} in Section \ref{s-2})
$$
\Hsup =\delta^{\beta_-}(1+\delta^\nu),  \mbox{ where $\nu<<1$.}
$$
Clearly $\usub_\eps-C\Hsup$ is also a local sub-solution of \eqref{original} in $\Omega_{\eps_0}$, where $C$ is an arbitrary positive number. Choose the value $C>0$ so large that $\usub_\eps-C\Hsup<\phi$  on $\Gamma_{\eps_0}$, that is
$$
L_\mu(\eps_0+\eps)-C\eps_0^{\beta_-}(1+\eps^\nu_0)<\min_{\Gamma_{\eps_0}} \phi.
$$
Because of the inequality $L_\mu(\eps_0+\eps)\leq L_\mu(\eps_0)$ the value $C=C(\eps_0)$ can be chosen independently of $\eps\in [0,\epsilon_0]$. With this fixed $C$ we now define the function
\begin{align}\label{sub-global}
\underline{\mathcal{U}}_\epsilon =
\begin{cases}
\max\{\usub_\eps-C\Hsup,\phi\} \tx{in} \Omega_{\eps_0},\\
\phi \tx{in} D_{\eps_0}.
\end{cases}
\end{align}
The function $\underline{\mathcal{U}}_\eps$ is a global sub-solution for all $\epsilon\in [0,\epsilon_0]$. Moreover
since $\Hsup =0$ on $\partial \Omega$ and $\usub_\eps$ is positive in $\Omega_{\eps_0}$, we have $\underline{\mathcal{U}}_\eps=\usub_\eps-C\Hsup$ near $\partial \Omega$.
Set $\omega_\eps:=\{ x\in \Omega_{\eps_0}: \usub_\eps-C\Hsup>\phi\}$ and note that $\omega_\eps\supset \omega_{\eps_0}$ for all $\eps\in [0,\eps_0]$, so that each $\omega_\eps$ contains a fixed neighbourhood of the boundary $\partial\Omega$. Thus
$$
\underline{\mathcal{U}}_\eps(x)=L_\mu(\delta(x)+\eps)-C(\eps_0)\delta^{\beta_-}(x)(1+\delta^\nu(x)) \mbox{ for } x\in\omega_{\eps_0} \mbox{ and for all } \epsilon \in (0,\eps_0).
$$
If $u$ is any solution of \eqref{original} which tends to infinity at the boundary
then the comparison principle of Lemma~\ref{comp2} implies that $u(x)\geq \underline{\mathcal{U}}_\epsilon(x)$ in $\Omega$ for all $\eps\in (0,\epsilon_0]$. Letting $\epsilon\to 0$ we get that
$u(x)\geq \underline{\mathcal{U}}_0(x)$ in $\Omega$ and in particular we find near the boundary that
$$
u(x)\geq \underline{\mathcal U}_0(x) = L_\mu(\delta(x))-C(\eps_0)\delta^{\beta_-}(x)(1+\delta^\nu(x))
$$
This together with \eqref{L_as} implies \eqref{limit}.
\qed

\begin{remark} If the domain $\Omega$ is small in the sense that its inradius $\rho_0$ satisfies
$$
\frac{\mu}{\rho_0^2}\geq e
$$
then $v=1$ is a global sub-solution. If $\mu\geq C_H(\Omega)$ it is not clear whether we can deduce from this fact that for large solutions $u$ the inequality $u>1$ holds.
\end{remark}

Proposition \ref{lbound} and \eqref{uAs} imply the following.
\begin{theorem} \label{bbe}  If $0<\mu< C_H(\Omega)$ then every large solution of \eqref{original} satisfies
\begin{align}\label{e-bbe}
\lim_{x \to \partial \Omega} \frac{u(x)}{\log \delta^{-2}(x)} =1.
\end{align}
\label{asymptotic}
\end{theorem}

\section{Uniqueness and existence of large solutions}\label{s-6}
\subsection{Uniqueness}

\begin{theorem}\label{uniqu}
Assume that $0<\mu<C_H(\Omega)$. Then \eqref{original} has at most one large  solution.
\end{theorem}
\proof
Suppose that \eqref{original} has two large solutions $U_1$ and $U_2$. If the domain is large they can become negative. In this case we add a sufficiently large negative multiple of the function $\phi\in H_0^1(\Omega)$ of Lemma~\ref{variational}
(recall that $\L_\mu \phi = -1$ and $\phi<0$ in $\Omega$) such that $w_i:=U_i -t\phi>1$ for $i=1,2$ and $t>0$ is taken sufficiently large. Then
$$
\L_\mu w_i=-\underbrace{a(x)}_{:=e^{t\phi(x)}}e^{w_i}+t \tx{in} \Omega, \quad w_i(x)\to \infty \tx{as} x\to \partial \Omega, \quad i=1,2.
$$
Define the function $\sigma(x)$ by $w_1(x)=\sigma(x) w_2(x)$.  Because of the asymptotic behavior of $U_1, U_2$ known from Theorem \ref{asymptotic} we have $\sigma(x)=1$ on $\partial \Omega$. Then
\begin{align*}
t&=\L_\mu w_1 +ae^{w_1}=-\sigma\Delta w_2 -w_2 \Delta \sigma -2\nabla \sigma\cdot \nabla w_2-\mu \delta^{-2} \sigma w_2  +ae^{\sigma w_2}\\
&=-w_2\Delta \sigma -2\nabla \sigma\cdot \nabla w_2  -\sigma ae^{w_2}+t\sigma +ae^{\sigma w_2}.
\end{align*}
Suppose that  $w_1>w_2$ (or equivalently $\sigma >1$) in a subset $\Omega'$ of $\Omega$. Since $w_1/w_2\to 1$ as $x$ approaches the boundary of $\Omega'$ we have $\sigma(x)=1$ on $\partial \Omega'$. Using our assumption $w_2>1$ we conclude that $e^{\sigma(x) w_2}>\sigma(x) e^{w_2}$ in $\Omega'$. Thus
$$
-w_2\Delta \sigma -2\nabla \sigma \cdot\nabla w_2 < t(1-\sigma)<0 \tx{in} \Omega',
$$
and by the maximum principle it follows that $\sigma \leq 1$ in $\Omega'$. This contradicts the fact that $w_1>w_2$ in $\Omega'$. Consequently we have $w_1\leq w_2$. Similarly we show that $w_2>w_1$ is impossible. This establishes the assertion.
\qed

\subsection{Existence}
\begin{theorem}\label{exist}
If $0<\mu<C_H(\Omega)$ then \eqref{original} has a large solution.
\end{theorem}
\proof
Let $\usup$ be a super-solution to \eqref{original} which blows up at $\partial\Omega$, as constructed in \eqref{e-usup}.
Let $\usub_m$ be a sub-solution to \eqref{original} defined in \eqref{sub-global} and chosen in such a way that $\usub_m=m$
on $\partial\Omega$ for $m\in \N$.
Let $\{M_n\}_{n\in \N}$ be a monotone increasing sequence of numbers satisfying
$$
\usub_m<M_n<\usup \tx{on} \Gamma_{1/n}.
$$
Let $u_{m,n}$ be the solution of the problem
$$
\L_\mu u_{m,n}+e^{u_{m,n}}=0\tx{in} D_{1/n},\quad u_{m,n}=M_n \tx{on} \partial D_{1/n}.
$$
Such a solution could be, e.g., constructed by minimizing the energy functional
$$J(u)=\frac{1}{2}\E_\mu(u)+\int_\Omega e^u\,dx,$$
which is coercive and weakly lower semicontinuous on the convex set
$$\mathcal{M}_n=\{u\in H^1(D_{1/n}),\;u=M_n\tx{on}\partial D_{1/n}\}.$$
>From the comparison principle of Lemma \ref{comp2} (i) it follows that
$$\usub_m\le u_{m,n}\le \usup\tx{in} D_{1/n}.$$
Thus, by standard compactness and diagonalization arguments we conclude that there exists a subsequence $\{u_{m,n(m)}\}_{m\in \N}$
which converges as $m\to\infty$  to a large solution $u$ of \eqref{original} in $\Omega$.
\qed

\section{Borderline potentials. Summary and open problems}\label{s-7}
By Theorem \ref{nonexistence}, no large solution of \eqref{original} exists if $\mu$ is negative.
This is due to the fact that the corresponding large sub-harmonics which interact with the nonlinear regime are too large near the boundary
and hence incompatible with the a priori bound constructed in Lemma \ref{Ob1}.
We are going to construct a maximal (in a certain sense) positive perturbation of $-\Delta$ of the form
$$\L_{\gamma(\delta)}:=-\Delta+\frac{\gamma(\delta)}{\delta^2},$$
where $\gamma(\delta)>0$, $\gamma(\delta)=o(1)$ as $\delta\to 0$, and
such that the semilinear problem
\begin{align}\label{gamma}
\L_{\gamma(\delta)} u+e^u=0\tx{in}\Omega
\end{align}
admits a large solution. Observe the different signs in the definition of $\L_{\gamma(\delta)}$ and $\L_\mu$.
Lemma \ref{Ob1} and the Phragmen--Lindel\"of alternative suggest that it is reasonable to look for a function $\gamma$
for which operator $\L_{\gamma(\delta)}$ admits large local sub-harmonics with the same or with a smaller order of magnitude as
the Keller--Osserman bound near $\partial \Omega$.
\smallskip

The asymptotic behavior given in \eqref{bobe} suggests to use
$$
h:=\Big(\log\frac{1}{\delta^{2}}\Big)^{m}, \quad m>0
$$
as a `prototype' family
of sub and super-harmonics in order to determine the borderline potential $\gamma(\delta)$.
By direct computations we have
$$
\L_{\gamma(\delta)}h=-\Delta h +\frac{\gamma(\delta)}{\delta^2} h = -h'\Delta \delta -h''(\delta)|\nabla \delta|^2 +\frac{\gamma(\delta)}{\delta^2}h,
$$
where $|\nabla \delta| =1$ and $\Delta \delta = -(N-1)\mathcal{H}_0 + o(\delta)$. Therefore
\begin{align*}
\L_{\gamma(\delta)}h=& \left\{\frac{2m}{\delta}\Big(\log\frac{1}{\delta^{2}}\Big)^{m-1}\Delta \delta
-\frac{4m(m-1)}{\delta^2}\Big(\log\frac{1}{\delta^{2}}\Big)^{m-2}|\nabla\delta|^2\right\}\\
&-\frac{2m}{\delta^2}\Big(\log\frac{1}{\delta^{2}}\Big)^{m-1}|\nabla\delta|^2+
\frac{\gamma(\delta)}{\delta^2}\Big(\log\frac{1}{\delta^{2}}\Big)^{m},
\end{align*}
where the expression in brackets is of lower order as $\delta\to 0$. Now we want to construct  $\gamma(\delta)$
such that $h$ is, depending on the value of $m$, either a sub or a super-harmonic.  Set
$$\gamma(\delta):=\beta\min\left\{\Big|\log\frac{1}{\delta^{2}}\Big|^{-1},1\right\},$$
for some $\beta>0$.
With such a choice of $\gamma$ we find that
$$
\L_{\gamma(\delta)}h=\frac{\beta-2m}{\delta^2}\Big(\log\frac{1}{\delta^{2}}\Big)^{m-1}\big(1+o(1)\big)
$$
in a small parallel strip $\Omega_\rho$. Therefore,
$$\Hsup:=\Big(\log\frac{1}{\delta^{2}}\Big)^{m}$$
is a local super-harmonic of $\L_{\gamma(\delta)}$ for all $0<m<\beta/2$. Otherwise, for $m>\beta/2$, $\Hsup$ is a local sub-harmonic of $\L_{\gamma(\delta)}$.
\smallskip

Further, a simple computation verifies that
$$\hsup=\delta^\alpha$$
is also a local super-harmonic of $\L_{\gamma(\delta)}$, for all $0\le\alpha<1$.
Thus, a Phragmen--Lindel\"of type argument similar to the one used in Theorem \ref{t-PL}, applied here to $\Hsup$ and $\hsup$
defined above,
shows that if $\hsub\ge 0$ is a local sub-harmonic of $\L_{\gamma(\delta)}$ then either
$$
(i) \quad  \limsup_{x\to \partial \Omega}\hsub(x)\Big(\log\frac{1}{\delta^{2}}\Big)^{-m}>0,
\quad\forall\, 0<m<\beta/2;
$$
or
$$(ii) \quad \quad \quad \hsub=0\tx{on}\partial\Omega.$$
In particular, every large solution of \eqref{gamma} must satisfy $(i)$ above.

Note that operator $\L_{\gamma(\delta)}$ is positive definite on $\Omega$, simply because $\gamma(x)>0$ in $\Omega$. As a consequence,
a comparison principle similar to Lemma \ref{comp2} (i) is valid for equation \eqref{gamma}.
Exactly the same arguments as in Lemma \ref{Ob1} imply that for large
 $A>0$ every solution $u$ of \eqref{gamma} satisfies
a Keller--Osserman type bound
\begin{align}\label{KO-gamma}
u(x)\leq \log \frac{A}{\delta^2(x)} \tx{in}\Omega.
\end{align}
Combining \eqref{KO-gamma} with the Phragmen--Lindel\"of bound $(i)$ which holds for any $m<\beta/2$, we immediately obtain a nonexistence result.
\begin{theorem} \label{non-gamma}
If $\beta>2$ then \eqref{gamma} does not have large solutions.
\end{theorem}

Next observe that  if $0<\beta <2$ then for $0<B<2-\beta$ the function
\begin{equation}
\usub=\log\frac{B}{\delta^2}
\label{l}
\end{equation}
is a local sub-solution of \eqref{gamma} with infinite boundary values. This local sub-solution can be extended to a global sub-solution in the same way as in \eqref{sub-global}. However, contrary to the construction in Proposition~\ref{lbound}, this time we cannot construct sub-solutions with everywhere finite and non-zero boundary values, cf. (i) in the conclusion from the Phragmen-Lindel\"of argument above.
\smallskip

In fact, we can prove the following existence and nonuniqueness result.
\begin{theorem}
If $0<\beta<2$ then \eqref{gamma} has a large solution $u$ such that
$$
\lim_{x \to \partial \Omega} \frac{u(x)}{\log \delta^{-2}(x)} =1,
$$
and for every $M>0$ there exists a large solution $v_M$ such that
$$
\lim_{x \to \partial \Omega} \frac{v_M(x)}{\big(\log \delta^{-2}(x)\big)^{\beta/2}} =M.
$$
\end{theorem}
\proof Recall that in Theorem \ref{exist} the existence was based on a family of sub-solutions with finite boundary values and a super-solution with infinite boundary value. Since such sub-solutions are no longer available in the present case, we sketch a different argument for the proof of the above existence result. For $k\in \N$ let $u_k$ be the large solution of
$$
\L_{\gamma(\delta)} u_k+e^{u_k}=0 \mbox{ in }D_{1/k}, \quad u_k = \infty \mbox{ on } \partial D_{1/k}.
$$
The sequence $u_k$ is monotonically decreasing, and if $\usub$ is the sub-solution from \eqref{l} extended to the whole of $\Omega$, then $u_k\geq \usub$ by the comparison principle. Therefore $u_k\to u$ as $k\to \infty$ locally uniformly in $\Omega$, where $u$ is a large solution of \eqref{gamma} in $\Omega$ with $u\geq \usub$. Hence $\lim_{x \to \partial \Omega} \frac{u(x)}{\log \delta^{-2}(x)} \geq 1$. Together with the Keller-Osserman upper bound from \eqref{KO-gamma} this establishes the first claim of the theorem.
\medskip

\newcommand{\vsup}{\overline{v}}
\newcommand{\vsub}{\underline{v}}
We now proceed to the construction of the large solution $v_M$. Let $M>0$ be any given number
and set
$$
\Hsub_{M,k}:=M\Big(\log\frac{1}{\delta^{2}}\Big)^{\beta/2}-k.$$
Then a straightforward computation yields for $\delta(x)$ small
\begin{multline*}
\L_{\gamma(\delta)}\Hsub_{M,k}+e^{\Hsub_{M,k}}\\
=\left\{M\beta(2-\beta)\delta^{-2}(\log(\delta^{-2}))^{\frac{\beta}{2}-2}(1+o(1))\right\}
-k\beta \delta^{-2}(\log(\delta^{-2}))^{-1}+e^{-k}e^{M(\log(\delta^{-2})^{\beta/2}}.
\end{multline*}
Since $\beta<2$ the expression in the parenthesis $\left\{\dots\right\}$ is of lower order as $\delta \to 0$. Let  $0<\eps<1$ and choose $\delta_0$ such that $M<(1-\epsilon)(\log(\delta_0^{-2}))^{1-\beta/2}$.
Then for all $x\in \Omega$ with $\delta(x)\leq \delta_0$ one finds
$$
\L_{\gamma(\delta)}\Hsub_{M,k}+e^{\Hsub_{M,k}} \leq -k\beta \delta^{-2}(\log(\delta^{-2}))^{-1}(1+o(1))+e^{-k}\delta^{-2(1-\eps)}\le 0
$$
provided $k>0$ is chosen sufficiently large. Hence $\Hsub_{M,k}$ is a local sub-solution.
Let $\phi\in H^1_0(\Omega)$ be defined as the solution of $\L_{\gamma(\delta)}\phi = -1$ (cf. Lemma \ref{variational} with $\mu$ replaced by $-\gamma(\delta)$). Similarly to \eqref{sub-global}, one can choose $k>0$ large enough so that $\vsub_M:=\max\{\Hsub_{M,k},\phi\}$
is a global sub-solution of \eqref{gamma} in $\Omega$.
\smallskip

To construct a super-solution, set
$$\Hsup_{M,K}:=M\Big(\log\frac{1}{\delta^{2}}\Big)^{\beta/2}+K,$$
which for large $K$ and $\delta(x)$ small is a local super-solution.
Let $A>0$ be as in Lemma \ref{Osserman}, so that $L_A(d(x))$ is a global super-solution of \eqref{gamma} in $\Omega$. Then one can choose $K>0$ so large that $\vsup_M:=\min\{L_A(d(x)),\Hsup_{M,K}\}$ is a global super-solution of \eqref{gamma} in $\Omega$, which coincides with $\Hsup_{M,K}$ near the boundary.

Since $\vsub_M<\vsup_M$ in $\Omega$,
a global large solution $v_M$ of \eqref{gamma} with the required asymptotic
can be constructed using a diagonalization procedure similar to the one
used in Theorem \ref{exist}.
We omit the details.
\qed

\subsection{Summary and open problems}

Our results are summarized as follows. Existence/nonexistence of large solutions for the problem
$$
-\Delta u -V(x)u +e^u = 0 \mbox{ in } \Omega, \quad u = 0 \mbox{ on } \partial\Omega
$$
can be read from the following table where we use the notation
$$\gamma_0=\min\left\{\Big|\log\frac{1}{\delta^{2}}\Big|^{-1},1\right\}.$$

\begin{center}
\begin{tabular}{|c|c|c|c|}
\hline
& $V(x)=\frac{\mu}{\delta^2}$ & $V(x)=\frac{-\beta \gamma_0(\delta)}{\delta^2}$ & $V(x)=\frac{\mu-\beta\gamma_0(\delta)}{\delta^2}$\\ \hline
\raisebox{-0.5em}{$\not\exists$} & \raisebox{-0.5em}{$\mu<0$} & \raisebox{-0.5em}{$\beta>2$} & $\mu<0$\\
& & & or $\mu = 0$, $\beta>2$ \\ \hline
\raisebox{-0.5em}{$\exists$} & \raisebox{-0.5em}{$0\leq \mu < C_H(\Omega)$} & \raisebox{-0.5em}{$0<\beta<2$} & $0<\mu<C_H(\Omega)$\\
& & & or $\mu=0$ and $0<\beta<2$. \\ \hline
critical & \raisebox{-0.5em}{no} & \raisebox{-0.5em}{$\beta=2$} & \raisebox{-0.5em}{$\mu=0$, $\beta=2$} \\
borderline &  & & \\ \hline
\end{tabular}
\end{center}
Except for $\mu=0$ the results of the last row in the above table were not proven in the present paper, but they can be obtained with little changes since for $\mu\not = 0$ the perturbation $\frac{\gamma_0(\delta)}{\delta^2}$ is of lower order than the dominant term $\frac{\mu}{\delta^2}$.
\smallskip

We finish our discussion with the following open questions:
\begin{enumerate}
\item Does  $\L_\mu u+ e^u=0$, $u\in H_0^1(\Omega)$ admit a solution for $\mu>C_H(\Omega)$ (see also Remark \ref{rem})?
\item Does \eqref{original} admit a large solution for $\mu\ge C_H(\Omega)$ ?
\item Does a solution of \eqref{original} exist with $u=\infty$ on $\Gamma_\infty$ and $u=0$ on $\Gamma_0$ where
$\Gamma_\infty \cup \Gamma_0= \partial \Omega$?
\item Does a large solution of \eqref{gamma} exist in the critical case $\beta=2$ ?
\end{enumerate}

\subsection*{Acknowledgements.}

Part of this work was supported by the Royal Society grant `Liouville theorems in nonlinear
elliptic equations and systems'. Part of the research was done while C.B. and V.M. were visiting the University of Karlsruhe (TH). The authors would like to thank this institution for its kind hospitality.

\end{document}